\documentstyle[12pt]{article}
\title{The Largest and Smallest Eigenvalues of Matrices and Some Hamiltonian Properties of Graphs}
\author {Rao Li \\                
         Dept. of Computer Science, Engineering, and Mathematics  \\
         University of South Carolina Aiken \\
	     Aiken, SC 29801 \\
	       USA  
         }
\date{Aug. 6, 2024}

\begin{document}
\maketitle
\begin{abstract}
Let $G = (V, E)$ be a graph. We define matrices $M(G; \alpha, \beta)$as $\alpha D + \beta A$,  where $\alpha$, $\beta$ are real numbers such that 
$(\alpha, \beta) \neq (0, 0)$ and 
$D$ and $A$ are the diagonal matrix and adjacency matrix of $G$, respectively. 
Using the largest and smallest eigenvalues of $M(G; \alpha, \beta)$ with $\alpha \geq \beta > 0$, we present sufficient conditions for the Hamiltonian and traceable graphs.
 \end{abstract} 
$$2010 \,\, Mathematics \,\, Subject \,\, Classification: \,\, 05C50, \, 05C45.$$
$$Keywords: the \,\, matrix, \,\, the \,\, largest \,\, eigenvalue, \,\, $$
$$\hspace*{24mm} Hamiltonian \,\,graph, \,\, traceable \,\,graph$$ 
\\

\noindent {\bf 1.  Introduction} \\

We consider only finite undirected graphs without loops or multiple edges.
Notation and terminology not defined here follow those in \cite{Bondy}.
For a graph $G = (V(G), E(G))$, we use $n$ and $e$ to denote its order $|V(G)|$ and size $|E(G)|$, respectively. The minimum degree and maximum degree
of $G$ are denoted by $\delta(G)$ and $\Delta(G)$, respectively. We use $N(u)$ to denote the set of all vertices which are adjacent to $u$ in $G$. 
A set of vertices in a graph $G$ is independent if the vertices in the set are pairwise nonadjacent. A maximum independent set in a graph $G$ is an independent set of largest possible size. 
 The independence number, denoted $\gamma(G)$, of a graph $G$ is the cardinality of a maximum independent set in $G$. 
For disjoint vertex subsets $X$ and $Y$ of $V(G)$,  
we define $E(X, Y)$ as $\{\, f : f = xy \in E, x \in X, y \in Y \,\}$.  
A cycle $C$  in a graph $G$ is called a Hamiltonian cycle of $G$ if $C$ contains all the vertices of $G$. 
A graph $G$ is called Hamiltonian if $G$ has a Hamiltonian cycle.
A path $P$  in a graph $G$ is called a Hamiltonian path of $G$ if $P$ contains all the vertices of $G$. 
A graph $G$ is called traceable if $G$ has a Hamiltonian path. \\

For a graph $G$, we define matrices $M(G; \alpha, \beta) := \alpha D + \beta A$, where $\alpha$, $\beta$ are real numbers such that 
$(\alpha, \beta) \neq (0, 0)$ and
$D$ and $A$ are the diagonal matrix and adjacency matrix of $G$, respectively. 
If $\alpha = 0$ and $\beta = 1$ (resp. $\alpha = 1$ and $\beta = 1$), then $M(G; \alpha, \beta)$ is the same as the adjacency matrix (resp. the signless Laplacian matrix) of $G$. Thus $M(G; \alpha, \beta)$ is a generalization of both adjacency matrix and  the signless Laplacian matrix of $G$. We use 
$\lambda_{\alpha, \beta; \, 1} \geq \lambda_{\alpha, \beta; \, 2} \geq \cdots \geq \lambda_{\alpha, \beta; \, n}$ to denote the eigenvalues of $M(G; \alpha, \beta)$. Since $M(G; \alpha, \beta)$ is symmetric, $\lambda_{\alpha, \beta; \, 1} \geq \lambda_{\alpha, \beta; \, 2} \geq \cdots \geq \lambda_{\alpha, \beta; \, n}$ are real numbers. 
Using the largest and smallest eigenvalues of $M(G; \alpha, \beta)$ with $\alpha \geq \beta > 0$, we in this note present sufficient conditions for the Hamiltonian and traceable graphs. The main results are as follows. \\

\noindent {\bf Theorem 1.} Let $G$ be $k$-connected ($k \geq 2$) a graph of $n \geq 3$ vertices and $e$ edges. Suppose $\alpha \geq \beta  > 0$.
Set $\lambda_1 := \lambda_{\alpha, \beta; \, 1}$ and $\lambda_n := \lambda_{\alpha, \beta; \, n}$.\\

\noindent [$1$] If 
$$\lambda_1 \leq (\alpha + \beta) \sqrt{\frac{(k + 1) \delta^2}{n}+ \frac{e^2}{n (n - k - 1)}},$$
then $G$ is Hamiltonian or $G$ is $K_{k, \, k + 1}$. \\

 \noindent [$2$] If 
 $$\lambda_n \geq (\alpha + \beta)\sqrt{\frac{(n - k - 1) \Delta^2}{n} + \frac{e^2}{n (k + 1)}},$$
then $G$ is Hamiltonian or $G$ is $K_{k, \, k + 1}$.\\

\noindent {\bf Theorem 2.} Let $G$ be $k$-connected ($k \geq 1$) a graph of $n \geq 9$ vertices and $e$ edges. Suppose $\alpha \geq \beta  > 0$.
Set $\lambda_1 := \lambda_{\alpha, \beta; \, 1}$ and $\lambda_n := \lambda_{\alpha, \beta; \, n}$.\\

\noindent [$1$] If 
$$\lambda_1 \leq (\alpha + \beta) \sqrt{\frac{(k + 2) \delta^2}{n}+ \frac{e^2}{n (n - k - 2)}},$$
then $G$ is Hamiltonian or $G$ is $K_{k, \, k + 2}$. \\

 \noindent [$2$] If 
 $$\lambda_n \geq (\alpha + \beta)\sqrt{\frac{(n - k - 2) \Delta^2}{n} + \frac{e^2}{n (k + 2)}},$$
then $G$ is Hamiltonian or $G$ is $K_{k, \, k + 2}$.\\

\noindent {\bf 2. Lemmas} \\

We will use the following results as our lemmas. The first two are from \cite{CE}. \\

\noindent {\bf Lemma $1$ } \cite{CE}. Let $G$ be a $k$-connected graph of order $n \geq 3$. If $\gamma \leq k$, then $G$ is Hamiltonian. \\

\noindent {\bf Lemma $2$} \cite{CE}. Let $G$ be a $k$-connected graph of order n. If $\gamma \leq k + 1$, then $G$ is traceable. \\

Lemma $3$ below is from \cite{M}.\\

\noindent {\bf Lemma $3$} \cite{M}. Let $G$ be a balanced bipartite graph of order $2n$ with bipartition ($A$, $B$). If
$d(x) + d(y) \geq n + 1$ for any $x \in A$ and any $y \in B$ with $xy \not \in E$, then $G$ is Hamiltonian. \\

Lemma $4$ below is from \cite{J}. \\ 

\noindent {\bf Lemma $4$} \cite{J}. Let $G$ be a $2$-connected bipartite graph with bipartition ($A$, $B$), where $|A| \geq |B|$. If
each vertex in $A$ has degree at least $s$ and each vertex in $B$ has degree at least $t$, then $G$ contains a cycle
of length at least $2 \min (|B|, s + t - 1, 2s - 2)$. \\

\noindent {\bf 3. Proofs} \\

\noindent{\bf Proof of Theorem 1.} Let $G$ be a $k$-connected ($k \geq 2$) graph with $n \geq 3$ vertices and $e$ edges.  Suppose $G$ is not Hamiltonian.
Then Lemma $1$ implies that $\gamma \geq k + 1$. Also, we have that $n \geq 2 \delta + 1 \geq 2 k + 1$ otherwise 
$\delta \geq k \geq n/2$
and $G$ is Hamiltonian.
 Let $I_1 := \{\, u_1, u_2, ..., u_{\gamma} \,\}$ be a maximum independent set in $G$. Then
  $I := \{\, u_1, u_2, ..., u_{k + 1} \,\}$ is an independent set in $G$.
  Thus
$$ \sum_{u \in I} d(u) = |E(I, V - I)| \leq \sum_{v \in V - I} d(v).$$ 
Since $\sum_{u \in I} d(u) +  \sum_{v \in V - I} d(v) = 2e$, we have that $$\sum_{u \in I} d(u) \leq e \leq \sum_{v \in V - I} d(v).$$
Let $V - I = \{\, v_1, v_2, ..., v_{n - (k + 1)}\,\}$. From Cauchy-Schwarz inequality, we have 
$$\sum_{r = 1}^{n - (k + 1)} 1^2 \,\, \sum_{r = 1}^{n - (k + 1)} d^2(v_r) \geq \left( \sum_{r = 1}^{n - (k + 1)} d(v_r)\right)^2 \geq e^2.$$
Thus 
$$\sum_{v \in V - I} d^2(v) \geq \frac{e^2}{n - k - 1}.$$
Therefore $$M := (k + 1) \delta^2 + \frac{e^2}{n - k - 1} \leq \sum_{u \in I} d^2(u) + \sum_{v \in V - I} d^2(v) = \sum_{v \in V} d^2(v).$$
with equality if and only if 
$d(u) = \delta$ for each $u \in I$, $\sum_{v \in V - I} d(v) = e$ (implying $\sum_{u \in I} d(u) = e$ and thereby $G$ is bipartite with partition sets of $I$ and $V - I$), 
and $\Delta = d(v)$ for each $v \in V - I$. \\

From Cauchy-Schwarz inequality again, we have 
$$\sum_{r = 1}^{k + 1} 1^2 \,\, \sum_{r = 1}^{k + 1} d^2(u_r) \leq \left( \sum_{r = 1}^{k + 1} d(u_r)\right)^2 \leq e^2.$$
Thus 
$$\sum_{u \in I} d^2(u) \leq \frac{e^2}{k + 1}.$$
Therefore $$N := \frac{e^2}{k + 1} + (n - k - 1) \Delta^2 \geq \sum_{u \in I} d^2(u) + \sum_{v \in V - I} d^2(v) = \sum_{v \in V} d^2(v).$$
with equality if and only if 
$d(v) = \Delta$ for each $v \in V - I$, $\sum_{u \in I} d(u) = e$ (implying $\sum_{v \in V - I} d(v) = e$ and thereby $G$ is bipartite with partition sets of $I$ and $V - I$), and $\delta = d(u)$ for each $u \in I$. \\

For any real row vector $X = (x_1, x_2, ..., x_n)$, we have 
$$XM(G; \alpha, \beta)X^T  = (\alpha - \beta) \sum_{i = 1}^{n} x_i^2 + \beta \sum_{uv \in E}(d(u) + d(v)^2 \geq 0,$$
where $X^T$ the transpose of $X$.   
Thus $M(G; \alpha, \beta)$ is positive semidefinite and therefore $\lambda_1 = \lambda_{\alpha, \beta; \, 1} \geq 
\lambda_{\alpha, \beta; \, 2} \geq \cdots \geq \lambda_{\alpha, \beta; \, n} = \lambda_n \geq 0$. Hence
$\lambda_1^2 = \lambda_{\alpha, \beta; \, 1}^2 \geq \lambda_{\alpha, \beta; \, 2}^2 \geq \cdots \geq \lambda_{\alpha, \beta; \, n}^2 
= \lambda_n^2 \geq 0$ are the eigenvalues of $M^2(G; \alpha, \beta)$.\\

Since $M^2(G; \alpha, \beta) = \alpha^2 D^2 + \alpha \beta D A + \alpha \beta A D + \beta^2 A^2$, the sum of all the entries in the $uth$ row of 
$M^2(G; \alpha, \beta)$ is equal to the sum of all the entries in the $uth$ rows of $\alpha^2 D^2$, $\alpha \beta D A$, $\alpha \beta A D$, and 
$\beta^2 A^2$, where $u$ is any vertex in $G$. Notice that
the sums of all the entries of the $uth$ rows of $D^2$, $D A$, $A D$, and $A^2$ are equal to $d^2(u)$, $d^2(u)$, $\sum_{v \in N(u)} d(v)$, and
$\sum_{v \in N(u)} d(v)$, respectively (See Page $805$ in \cite{LP}). Hence the sum of all the entries in the $uth$ row, denoted $RS(u)$, in $M^2(G; \alpha, \beta)$ is
$$\alpha (\alpha + \beta) d^2(u) + \beta (\alpha + \beta)\sum_{v \in N(u)} d(v).$$
Let $Y = (1, 1, ..., 1)$ be an $n$-dimensional row vector. Applying Rayleigh-Ritz theorem to $M^2(G; \alpha, \beta)$, we have 
$$\lambda_1^2 \geq \frac{YM(G; \alpha, \beta)Y^T}{YY^T} \geq \lambda_n^2.$$
Notice that $$YM(G; \alpha, \beta)Y^T = \sum_{u \in V} RS(u)= \alpha (\alpha + \beta) \sum_{u \in V}d^2(u) + \beta (\alpha + \beta)\sum_{u \in V} \sum_{v \in N(u)} d(v)$$
$$= \alpha (\alpha + \beta) \sum_{u \in V}d^2(u) + \beta (\alpha + \beta) \sum_{u \in V}d^2(u) = (\alpha + \beta)^2 \sum_{u \in V}d^2(u).$$ 
So 
$$\lambda_1^2 \geq (\alpha + \beta)^2 \,\, \frac{\sum_{u \in V}d^2(u)}{n} \geq \lambda_n^2.$$ 
\noindent [$1$] From the condition in this case, we have 
$$(\alpha + \beta)^2\left(\frac{(k + 1) \delta^2}{n} + \frac{e^2}{n(n - k - 1)}\right)$$
$$\geq \lambda_1^2 \geq (\alpha + \beta)^2 \frac{\sum_{u \in V} d^2(u)}{n} \geq (\alpha + \beta)^2 \frac{M}{n} $$
$$= (\alpha + \beta)^2\left(\frac{(k + 1) \delta^2}{n} + \frac{e^2}{n(n - k -1)}\right).$$
Thus each inequality above becomes an equality. Therefore $d(u) = \delta$ for each $u \in I$, $\sum_{v \in V - I} d(v) = e$ (implying $\sum_{u \in I} d(u) = e$ and thereby $G$ is bipartite with partition sets of $I$ and $V - I$), and and $\Delta = d(v)$ for each $v \in V - I$.  Hence
$$(k + 1) \delta = |E(I, V - I)| = \Delta (n - k - 1) \geq \delta (n - k - 1).$$ 
So $2k + 2 \geq n \geq 2k + 1$.  If $n = 2k + 2$, then $\delta = \Delta$. Lemma $3$ implies $G$ is Hamiltonian, a contradiction. 
If $n = 2k + 1$, then $G$ is $K_{k, \, k + 1}$. \\

This completes the proof of [$1$] in Theorem [$1$]. \\

\noindent [$2$] From the condition in this case, we have 
$$(\alpha + \beta)^2\left(\frac{(n - k - 1) \Delta^2}{n} + \frac{e^2}{n(k + 1)}\right)$$
$$\leq \lambda_n^2 \leq (\alpha + \beta)^2 \,\, \frac{\sum_{u \in V} d^2(u)}{n} \leq (\alpha + \beta)^2 \frac{N}{n} $$
$$= (\alpha + \beta)^2\left(\frac{(n - k - 1) \Delta^2}{n} + \frac{e^2}{n(k + 1)}\right).$$
Thus each inequality above becomes an equality. Therefore $d(v) = \Delta$ for each $v \in V - I$, $\sum_{u \in I} d(u) = e$ (implying $\sum_{v \in V - I} d(v) = e$ and thereby $G$ is bipartite with partition sets of $I$ and $V - I$), and $\delta = d(u)$ for each $u \in I$. 
Hence
$$(k + 1) \delta = |E(I, V - I)| = \Delta (n - k - 1) \geq \delta (n - k - 1).$$ 
So $2k + 2 \geq n \geq 2k + 1$.  If $n = 2k + 2$, then $\delta = \Delta$. Lemma $3$ implies $G$ is Hamiltonian, a contradiction. 
If $n = 2k + 1$, then $G$ is $K_{k, \, k + 1}$. \\

This completes the proof of [$2$] in Theorem [$1$]. \\

By using Lemma $2$, Lemma $3$, Lemma $4$, and the ideas in the proof of Theorem $1$ above, we can prove Theorem $2$.
 The details of the proof of Theorem $2$ are skipped. \\ 
 
\noindent {\bf Remark. }  From the proof of Theorem $1$, we have the following corollary. \\

\noindent {\bf Corollary. } Let $G$ be a graph of order $n$ vertices and $e \geq 1$ edges. Suppose $\alpha \geq \beta  > 0$ and $I$ is any  
independence set of $G$ with $|I| = \gamma$. Set $\lambda_1 := \lambda_{\alpha, \beta; \, 1}$ and $\lambda_n := \lambda_{\alpha, \beta; \, n}$.
Then\\

\noindent [$1$]  
$$\lambda_1 \geq (\alpha + \beta) \sqrt{\frac{\gamma \delta^2}{n}+ \frac{e^2}{n (n - \gamma)}}.$$

 \noindent [$2$] 
 $$\lambda_n \leq (\alpha + \beta)\sqrt{\frac{(n - \gamma) \Delta^2}{n} + \frac{e^2}{n \gamma}}.$$

\noindent {\bf Acknowledgment} \\

The author would like to thank the referees for their suggestions or comments leading to the improvements of the initial version of the manuscript.

\end{document}